\documentclass[conference]{worldcomp}

\usepackage[hmargin=.75in,vmargin=1in]{geometry}
\usepackage[american]{babel}
\usepackage[T1]{fontenc}
\usepackage{times}
\usepackage{caption}
\usepackage{epsfig}
\def\by{{\mathbf y}}
\def\bx{{\mathbf x}}
\def\bv{{\mathbf v}}
\def\bu{{\mathbf u}}
\def\bc{{\mathbf c}}
\def\be{{\mathbf e}}
\def\cS{{\mathcal S}}
\def\diag{\text{{diag}}}
\def\tr{\text{{tr}}}
\usepackage{textcomp}
\usepackage{epsfig,graphicx}
\usepackage{xcolor}
\usepackage{amsfonts,amsmath,amssymb}
\usepackage{fixltx2e} % Fixing numbering problem when using figure/table* 
\usepackage{booktabs}

\title{\bf Greedy Approach for Subspace Clustering from Corrupted and Incomplete Data}           %%%% Replace with your title.

\author{
{\bfseries A. Petukhov$^1$, I. Kozlov$^2$}\\
$^1$Contact Author, Department of Mathematics, University of Georgia, Athens, GA 30602, USA,\\ petukhov@math.uga.edu\\
$^2$Algosoft Tech USA, Bishop, GA, USA, inna@algosoft-tech.com
}

\begin{document}

\maketitle                        %%%% To set Title and Author names.

\begin{abstract}%%%% Replace with your abstract.
We describe the Greedy Sparse Subspace Clustering (GSSC)  algorithm providing an efficient method 
for clustering data belonging to a few low-dimensional linear or affine subspaces  from incomplete corrupted and noisy data. 
\par
We provide numerical evidences that, even in the simplest implementation, the greedy approach  
increases the subspace clustering capability of the existing state-of-the art SSC algorithm significantly.
\end{abstract}

\vspace{1em}
\noindent\textbf{Keywords:}
 {\small  Subspace Clustering, Sparse Representations, Greedy Algorithm, Law-Rank Matrix Completion, Compressed Sensing} %%%% Replace with your keywords

%%%%%%%%%%%%%%%%%%%%%%%%%%%%%%%%%%%%%%%%%%%%%%%%%%%%%%%%%%%%

\section{Introduction}
We consider a greedy strategy  based algorithm for preprocessing on vector database 
necessary for subspace clustering. The problem of subspace clustering consists in classification of 
the vector data belonging to a few linear or affine low-dimensional subspaces of the high dimensional 
ambient space when neither subspaces or even their dimensions are not known, i.e., they have to be identified  from the same database. 
\par
The algorithm considered below does not require clean input data. Vice versa, we assume that those data are corrupted with sparse errors, random noise distributed over all vector entries and quite significant part of data is missed. 
\par
The property of errors to constitute a sparse set means that some (but not all) vector entries
are corrupted, i.e., those values are randomly replaced. The locations (indices) of corrupted entries are unknown. 
\par
The noise is randomly introduced in each entry. Its magnitude is usually much less than the data magnitude.
\par
In information theory, the missed samples are called erasures. They have two main features.
First, the data values in erasures does not have practical importance. The most natural way is to think about erasures as about lost data. Second, the coordinates of erasures are known. 
\par
Using more formal definition, we have $N$ vectors $\{\by_l\}_{l=1}^N$ in $K$ linear or affine subspaces $\{\cS^l\}_{l=1}^K$ with the dimensions $\{d_l\}_{l=1}^K$ of the $D$-dimensional Euclidean space $\mathbb R^D$.  We do not assume that those spaces do not have non-trivial intersections. However, we do assume that any one of those spaces is not a subspace of other one.
At the same time, the situation when one subspace is a subspace of the direct sum of two other subspaces is allowed. Such settings  inspire the hope that when $N$ is large enough and the points are randomly and independently distributed on those planes some sophisticated algorithm can identify
those planes and classify the belongingness of each point to the found subspaces. Of course, some of points may belong to the intersection of two or more subspaces, then such point maybe assigned to
one of those space or to all of them. With the probability 1 the points belong to only one of subspaces.
Then the problem consists in finding a permutation matrix $\Gamma$ such that 
$$
[Y_1, \dots, Y_K]=Y\Gamma,
$$
where $Y\in \mathbb R^{D\times N}$ is an input matrix whose columns are the given points in an arbitrary random order, whereas  in $[Y_1, \dots, Y_K]$ is rearrangement of the matrix $Y$ in the accordance with the affiliation of the vector with the subspaces $\cS^k$. 
\par
The problem of finding clusters $\{Y_k\}$ is usually by means of finding the clusters in the graph whose edges (to be more precise the weights of edges) characterize the interconnection between pairs of vertexes. In our  case, the popular method of clustering consists in making the points to play role of the vertexes, while the weights are  set from the coefficients of decomposition of the vectors through other vectors from the same space $\cS^k$. This idea looks as vicious circle. We are trying to identify the space $\cS^k$ accommodating the vector $\by_i$, using its linear decomposition  in the remaining vectors of $\cS^k$. However, the situation is not hopeless at all. In   \cite{EV}, the excellent suggestion was formulated. Obviously, the decomposition problem formulated above is reduced to solving the non-convex problem 
\begin{equation}
\label{CS}
\|\bc_i\|_0\to\min,\text{ subject to }\by_i=Y\bc_i,\, c_{i,i}=0, \,i=1,\dots,n;
\end{equation}
where $\|\bx\|_0$  is the Hamming weight of the vector $\bx$, $\|\bx\|_0=\#\{x_j\ne0\}$. The problem of finding the sparsest solutions to (\ref{CS}), so-called, {\it Compressed Sensing} or {\it Compressive Sampling} underwent thorough study originated in \cite{CT}, \cite{D}, \cite{RV} and continued in hundreds of theoretical and applied papers.
\par
In the ideal world of perfect computational precision and unlimited computational power, with probability one, the solutions 
of (\ref{CS}) would point out the elements of the appropriate $\cS^l$ by their non-zero decomposition coefficients. Provided that no vectors of wrong subspaces participate in decomposition of each column of $Y$, the matrix $C$ whose columns are $\bc_i$ allows perfectly reconstruct the structure of the subspaces in polynomial time. There are two obstacles on that way. 
\par
First, the precision of the input matrix $Y$ is not perfect. So the decompositions may pick up wrong vectors even if we are able to solve problem (\ref{CS}).  In this case, the problem of subspace clustering is considered for the similarity graph defined by the symmetric matrix $W:=|C|+|C^T|$. While, generally speaking, this problem has non-polynomial complexity, there exist practical algorithms allowing right clustering when the number 
of "false" interconnections of elements from different subspaces are not very dense and not very intensive. Following \cite{EV}, we  use some modification of  the spectral clustering algorithm from \cite{NWJ} which is specified in \cite{L} as "graph's random walk Laplacian".
\par
The second obstacle consists in non-polynomial complexity of the problem (\ref{CS}). The elegant solution allowing to overcome this obstacle is replacement of  non-convex problem  (\ref{CS}) with the convex problem
\begin{equation}
\label{CS1}
\|\bc_i\|_1\to\min,\text{ subject to }\by_i=Y\bc_i,\, c_{i,i}=0, \,i=1,\dots,n;
\end{equation}
or
\begin{equation}
\label{CSM}
\|C\|_1\to\min,\text{ subject to }Y=YC,\, \text{ diag}(C)=0;
\end{equation}
in the matrix form. It follows from the fundamental results from  \cite{CT}, \cite{D}, \cite{RV} that for matrices $Y$ with some
reasonable restrictions on $Y$ and for not very large Hamming weight of the ideal sparse solution, it can be uniquely found by solving convex problem (\ref{CS1}). There are other more efficient $\ell^1$-based methods for finding sparse solutions (e.g., see \cite{CWB},  \cite{PK2}, \cite{KP1}). 
\par
It should be mentioned, that the matrices $Y$ in practical problems may be very far from the requirements  for the uniqueness of solutions. At the same time, the uniqueness of the solutions are not necessary in our settings. We just wish to have the maximum of separation between indices of 
the matrix $W$ corresponding to different subspaces. 
\par
In the case of successful clustering, the results for each $\cS^k$ may be used for further processing like
data noise removal, error correction, and so on. Such procedures become significantly more efficient when applied to low-rank submatrices of $Y$ corresponding to one subspace.
\par
Thus, in applications, the problems involving subspace clustering can be split into 3 stages: 1) preprocessing; 2) search for clusters in the graphs; 3) processing on clusters.
In this paper,  we develop a first stage algorithm helping to perform the second stage much more efficiently then the state-of-the-art algorithms. We do not discuss any aspects of improvements of stage 2.  We just 
take one of such algorithms, specifically the spectral clustering, and use it for comparison of the influence of our and competing preprocessing algorithms on the efficiency of clustering.
\par
As for stage 3, its content depends on an applied problem requesting subspace clustering. One of typical possible goal of the third stage can be data recovery from incomplete and corrupted measurements. Sometimes this problem  is called "Netflix problem". We will discuss below how the same problems of incompleteness and corruption can be solved within clustering preprocessing. However, for low-rank matrices it can be solved more efficiently. Among many existing algorithms we mention the most recent papers \cite{CR}, \cite{CJSC}
\cite{LCM}, \cite{WYGSM}, \cite{YYO}, \cite{PK4} providing the best results for input having both erasures and errors.
\par
In Section \ref{disc}, we discuss the problem formal settings. In Section \ref{alg}, the existing SSC algorithm and our modification will be given. The results of numerical experiments showing the consistency of the proposed approach will be given in Section \ref{exper}.
%%%%%%%%%%%%%%%%%%%%%%%%%%%%%%%%%%%%%%%%%%%%%
\par
 \section{Problem Settings}
 \label{disc}
We use Sparse Subspace Clustering algorithm (SSC) from \cite{EV} as a basic algorithm 
for our  modification based on a greedy approach. Therefore, significant part of reasoning we give in this
section can be found in  \cite{EV} and in the earlier paper \cite{EV0}. Very similar ideas of subspace self-representation  for subspace clustering were used also in \cite{RTVM}. However, the error resilience mechanism in that paper was used under assumption that there are enough uncorrupted data vector, whereas,  this assumption is not required in \cite{EV}.

\par
Optimization CS problems \ref{CS1} and  \ref{CSM} assume that the data are clean, i.e., they have no noise and errors. Considering the problem within the standard CS framework, the problem  \ref{CS} can be reformulated as finding the sparsest vectors $\bc$ (decomposition coefficients) and $\mathbf e$ satisfying the system of linear equations $\by=A\bc+\mathbf e$.
\par
It was mentioned in \cite{WYGSM} that the last system can be re-written as 
\begin{equation}
\label{errorcor}
\by=[A\, I]\left[\begin{array}{c} \bc \\ \mathbf e\end{array}\right],
\end{equation}
where $I$ is the identity matrix. Therefore, the problem of sparse reconstruction and error correction can be solved simultaneously with CS methods. In \cite{PK3}, we designed an algorithm efficiently finding sparse solutions to  (\ref{errorcor}).
\par
Unfortunately, the subspace clustering  cannot use this strategy straightforwardly because not only
"measurements" $\by_i$ are corrupted in (\ref{CS}) but "measuring matrix" $A=Y$ also can be corrupted.
It should be mentioned that if the error probability  is so  low that there exist uncorrupted columns of $Y$ constituting bases for all subspaces $\{\cS^k\}$, the method from \cite{PK3} can solve
the problem of sparse representation with simultaneous error correction. In what follows, the considered algorithm will admit a significantly higher error rate. In particular, all columns of $Y$ may be corrupted.
\par
Following \cite{EV}, we introduce two (unknown for the algorithm) $D\times N$  matrices $\hat E$ and $\hat Z$. The matrix $\hat E$ contains a sparse (i.e., $\#\{E_{ij}\ne 0\}<DN$) set of  errors with relatively large magnitudes.  The matrix $\hat Z$ defines the noise having relatively low magnitude but distributed over all entries of $\hat Z$.
Thus, the clean data are representable as $Y-\hat E -\hat Z$. Therefore, when the data are corrupted with sparse errors and noise, the equation $Y=YC$ has to be replaces by 
\begin{equation}
\label{corr}
Y=YC +\hat E(I-C)+\hat Z(I-C).
\end{equation} 
The authors of \cite{EV} applied a reasonable simplification of the problem by replacing 2 last terms 
of (\ref{corr}) with some (unknown) sparse matrix $E:=\hat E(I-C)$ and the matrix with the deformed noise $Z=\hat Z(I-C)$. Provided that sparse $C$ exists, the matrix $E$ still has to be sparse. This transformation  leads to some simplification of the optimization procedure. This is admissible simplification since, generally speaking, we do not need to correct and denoise the input data $Y$. Our only goal is  to find the sparse matrix $C$. Therefore, we do not need matrices $\hat E$ and $\hat Z$. 
While, as we mentioned above, the error correction procedure can be applied after subspace clustering, original setting (\ref{corr}) with the genuine values of the errors and noise within subspace clustering still makes sense and deserves to become a topic for future research.
\par
Taking into account modifications from last paragraph, the authors of \cite{EV} formulate constrained optimization problem 
\begin{equation}
\label{constr1}
\begin{gathered}
\min \|C\|_1+\lambda_e\|E\|_1+\frac{\lambda_z}{2}\|Z\|_{F}^2
\\
\text{s.t. }Y=YC+E+Z, \,\text{diag}(C)=0,
\end{gathered}
\end{equation}
where $\|\cdot\|_F$ is the Frobenius norm of a matrix. If the clustering into affine subspaces is required, the additional constrain $C^T\mathbf 1=\mathbf 1$ is added.
\par
On the next step, using the representation $Z=Y-YC-E$ and introducing an auxiliary matrix 
$A\in\mathbb R^{N\times N}$, constrained optimization problem (\ref{constr1}) is transformed into
\begin{equation}
\label{constr2}
\begin{gathered}
\min \|C\|_1+\lambda_e\|E\|_1+\frac{\lambda_z}{2}\|Y-YA-E\|_{F}^2
\\
\text{s.t. } A^T\mathbf 1=\mathbf 1,\, A=C-\text{diag}(C),
\end{gathered}
\end{equation}
Optimization problems  (\ref{constr1})  and (\ref{constr2})  are equivalent.
Indeed, obviously, at the point of extremum of (\ref{constr2}), diag$(C)=0$. Hence, $A=C$.
\par
At last, the quadratic penalty functions with the weight $\rho/2$ corresponding to constrains are added to
the functional in  (\ref{constr2}) and the Lagrangian functional is composed. The final Lagrangian functional is as follows
\begin{equation}
\label{lagr}
\begin{gathered}
\mathcal L(C,A,E, \boldsymbol \delta, \Delta)=\min \|C\|_1\\+\lambda_e\|E\|_1 +\frac{\lambda_z}{2}\|Y-YA-E\|_{F}^2
\\
+\frac{\rho}{2}\|A^T\mathbf 1-\mathbf 1\|_2^2+\frac{\rho}{2}\|A-C+\text{diag}(C))\|_F^2
\\
+\boldsymbol\delta^T(A^T\mathbf 1-\mathbf 1)+\tr(\Delta^T(A-C+\text{diag}(C))),
\end{gathered}
\end{equation}
where the vector $\boldsymbol \delta$ and the matrix $\Delta$ are Lagrangian coefficients.
Obviously, since the penalty functions are formed from the constrains, they do not change the point and value of the minimum. 
%%%%%%%%%%%%%%%%%%%%%%%%%%%%%%%%%%%%%%%%
\section{Algorithm}
\label{alg}
For minimization of functional (\ref{lagr}) an Alternating Direction Method  of Multipliers (ADMM,  \cite{BPCPE}) is used. In \cite{EV}, this is a crucial part of the entire algorithm which is called the Sparse Subspace Clustering algorithm. 
\par
The parameters $\lambda_e$ and $\lambda_z$ in (\ref{lagr}) are selected in advance. They define the compromise
between good approximation of $Y$ with $YC$ and the high sparsity of $C$. The general rule is 
to set the larger values of the parameters for the less level of the noise or errors. In \cite{EV}, the selection of the parameters 
by formulas
\begin{equation}
\label{pars}
	\lambda_e=\alpha_e/\mu_e,\qquad \lambda_z=\alpha_z/\mu_z,
\end{equation}
where $\alpha_e,\, \alpha_z>1$ and 
$$
\mu_e:=\min_i\max_{j\ne i}\|\by_j\|_1,\quad\mu_z:=\min_i\max_{j\ne i}|\by_i^T\by_j|,
$$
is recommended.
\par
The initial parameter $\rho=\rho^0$ is set in advance. It is updated as $\rho:=\rho^{k+1}=\rho^k\mu$ with iterations of SSC algorithm. We notice that, adding the penalty terms, we do not change the problem. It still has the same minimum. However,
the appropriate selection of $\mu$ and $\rho^0$ accelerates the algorithm convergence significantly.
\par
  Each iteration of the  algorithm is based on consecutive optimization with respect to each of the unknown values $A$, $C$, $E$, $\boldsymbol\delta$, $\Delta$ which are initialized by zeros before the algorithm starts. 
\par
Due to appropriate form of functional (\ref{lagr}), optimization of each value is simple and computationally efficient. The five formulas for updating the unknown values are discussed below.
\par
The matrix $A^{k+1}$ is a solution of the system of linear equations
\begin{multline}
\label{a}
(\lambda_zY^TY+\rho^kI+\rho^k\mathbf 1\mathbf 1^T)A^{k+1}=\lambda_zY^T(Y-E^k)\\
+\rho^k(\mathbf 1\mathbf 1^T+C^k)-\mathbf 1\boldsymbol\delta^{kT}-\Delta^k
\end{multline}
\par
When the data are located on linear subspaces the terms $\mathbf 1\mathbf 1^T$ and
$\mathbf 1\boldsymbol\delta^{kT}$ may be removed (set to 0) from (\ref{a}).
\par
While the system (\ref{a}) has matrices of size $N\times N$, due to its special form,
 the complexity of the algorithm
for the inverse matrix is $O(D^2N)$ that is much lower than $O(N^3)$, provided that $D\ll N$.
Unfortunately, the matrix $\Delta$ may have a full rank. Therefore, the computational cost of
its product with the inverse matrix is $O(DN^2)$, i.e., not so impressive as for the matrix inversion.
\par
We will need the following notation
$$
 S_{\epsilon}[x]:=\left\{\begin{array}{ll} x-\epsilon, &x>\epsilon,\\x+\epsilon, &x<-\epsilon,\\0, &\text{otherwise};\end{array}\right.
$$
where $x$ can be either a number or a vector or a matrix. The operator $S_{\epsilon}[\cdot]$ is called the shrinkage operator. 
\par
Let 
$$
J:=S_{\frac{1}{\rho}}[A^{k+1}+\Delta^k/\rho].
$$
Then the matrix $C^{k+1}$ is defined by the formula
\begin{equation}
\label{c}
C^{k+1}:=J-\diag(J).
\end{equation}
The remaining values $E^{k+1}$, $\boldsymbol\delta$, and $\Delta$ are computed
as
\begin{equation}
\label{e}
E^{k+1}:=S_{\frac{\lambda_e}{\lambda_z}}[Y-YA^{k+1}].
\end{equation}

\begin{equation}
\label{d}
\boldsymbol\delta^{k+1}:=\boldsymbol\delta^{k}+\rho^k(A^{k+1}\mathbf 1-\mathbf 1),
\end{equation}

\begin{equation}
\label{D}
\Delta^{k+1}:=
\Delta^{k}+\rho^k(A^{k+1}-C^{k+1}).
\end{equation}
\par
The algorithm goes to the next iteration if one of conditions
$$
\|A^{k+1}\mathbf 1-\mathbf 1\|_\infty<\epsilon,\quad
\|A^{k+1}-C^{k+1}\|_\infty<\epsilon,
$$
$$
\|A^{k+1}-A_k\|_\infty<\epsilon,\quad
\|E^{k+1}-E^k\|_\infty<\epsilon
$$
fails, where $\epsilon$ is the given error tolerance.
\par
In the form shown above SSC algorithm gives the state-of-the-art benchmarks for subspace clustering 
problems. 
\par
Our suggestion is to attract ideas of greedy algorithms to increase the capability of SSC algorithm in
subspace clustering.  Greedy algorithms are very popular in non-linear approximation (especially in redundant systems) when the global optimization is replaced with iterative selection of the most probable candidates
from the point of view their prospective contribution into approximation. The procedure is repeated 
with selection of new entries, considering the previously selected entries as reliable  with guaranteed participation in approximation. The most typical case is Orthogonal Greedy Algorithm consisting in selection of the approximating entries having the biggest inner products with the current approximation residual and follow-up orthogonal projection of the approximated object onto the span 
of selected entries. 
\par
In many cases, OGA allows to find the sparsest representations if they exist. In \cite{KP1} and \cite{PK2}, we applied greedy idea in combination with the reweighted $\ell^1$-minimization to CS problem of finding the sparsest solutions of underdetermined system.   We used the existing $\ell^1$-minimization scheme with the the opportunity to reweight entries. When the basic algorithm fails, the greedy blocks picks the biggest (the most reliable) entries in the decomposition whose magnitudes
are higher than some threshold. Those entries are considered as reliable. Therefore they get the less
weight in the $\ell^1$ norm while other entries are competing on next iterations for the right to be picked up. 
\par 
The similar idea was employed in our recent paper \cite{PK4}, where the greedy approach was applied 
to the algorithm  for completion of low-rank matrices from incomplete highly corrupted samples from \cite{LCM} based  on  Augmented Lagrange Multipliers method. The simple greedy modification of
the matrix completion algorithm from \cite{LCM} gave the boost in the algorithm restoration capability.
\par
Now we discuss details how the greedy approach can be incorporated in (to be more precise {\it over}) SSC algorithm.
First of all we introduce  a non-negative matrix $\Lambda\in\mathbb R^{D\times N}$ whose entries reflect our knowledge about the entries of error matrix $E$. Let us think that the regular entries
with no (say, side)  information have values $1$, whereas  entries with coordinates of presumptive errors are  set to small value or to 0.
\par
Let us consider the mechanism of the influence of the parameter $\lambda_e$ on the  output matrix $C$.  $\lambda_e$ sets
the balance between the higher level of the sparsity of $C$ with the more populated error matrix
$E$ and less sparse $C$ but less populated matrix $E$. 
Setting too small $\lambda_e$  makes too many "errors" and very sparse $C$. However,  probably, this is not what we want. This would mean that sake of $C$ sparsity we introduced too large distortion into the input data $Y$. At the same time, if we know for sure or almost for sure that some 
entry of $Y$ with indices $(i, j)$  is corrupted, we  loose nothing by assigning to this element an individual small weight in functional (\ref{lagr}). This weight can be much less than $\lambda_e$ or even equal to $0$. Thus, we have to replace the term $\|E\|_1$ with $\|\Lambda\odot E\|_1 $, where operation $\odot$ means entrywise product of two matrices. In practice,
this means that in formula (\ref{e}) we will apply different shrinkage threshold for different indices.
Generally speaking, it makes sense to use all range of non-negative real numbers to reflect our knowledge about $E$. Say, highly reliable entries have to be protected from distortion by the weight greater than 1 in $\Lambda$. However, in this paper we restrict ourself with two-level entries:
either 1 (no knowlege) or $10^{-4}$ (suspicious to be an error). 
\par
 When no a priori knowledge is available, we set all entries of $\Lambda$ to 1. However, 
in information theory there is a special form of corrupted information which is called {\it erasure}. The losses of network packets is the most typical reason of  erasures. Another example of erasure is given by occlusions in video models when moving object temporary overlap the background or each other. 
Erasures represent rather missing than corrupted information.  Erased entries like entries with errors have to be restored or at least taken into account.   The only difference of erasures with errors
is a priory knowledge of their locations. This additional information allows to reconstruct the values of erasures  in more efficient way than errors. So the entries with erasures  have to be marked (say, by setting the corresponding entries of $\Lambda$ to 0) before the algorithm starts. 
\par
The entries which are suspicious to be  errors  dynamically extend "the  map of marked errors" in $\Lambda$ after each iteration of the greedy algorithm.
\par
Thus, in the Greedy Sparse Subspace Clustering (GSSC) we organize an external loop over sparsification part of SSC algorithm with an additional  input matrix $\Lambda$ which is initialized  with ones at regular entries and with some small non-negative number $\kappa\ge0$ at the places of erased entries. 
\par
One iteration of our greedy algorithm consists in running the modified version of SSC and  
$\Lambda$ and $Y$ updates.
Our modification of SSC consists of two parts. First, we  take into account the matrix $\Lambda$ while computing $E^{k+1}$ by formula (\ref{e}). Second,  we  update $\rho$ on each iteration of greedy algorithm.
\par
After the first iteration, the estimated matrix  $E$ is used to set the threshold 
$$
T^1=\max(\alpha_1\|Y-E\|_\infty, \,\alpha_2\max_{1\le j \le N}\text{\,median\,} |\by_j|), 
$$
where $ 0<\alpha_1,\alpha_2<1$.
 We use the median estimate to avoid too low threshold leading to  large error map when there is no error in the data.
\par
Starting from the second greedy iteration, we just update threshold $T^{n+1}=\beta T^n$, $0<\beta<1$.
The current value $T^{n+1}$ is used for the extension of the error map by formula
$$
\Lambda^{n+1}_{ij}=
\left\{
\begin{array}{ll}
\Lambda^n_{ij}, &|E_{ij}^{n+1}|<T^{n+1},
\\
\kappa, &|E_{ij}^{n+1}|\ge T^{n+1}\end{array},
\right.
$$
where  $E^{n+1}$ is the error matrix obtained on the previous iteration.
The updated $\Lambda$ is used for next iteration. 
\par
In addition, on each iteration of greedy algorithm, we update the input matrix of SSC algorithm by the formula $Y_{ij}^{n+1}=Y_{ij}^n-E_{ij}^n$ for the pairs $(i,j)$ marked in the current $\Lambda^n$ as errors/erasures.
While $E^n$ is not a genuine matrix of errors,  this is not serious drawback for the original SSC algorithm. However,
for GSSC this may lead to unjustified and very undesirable extension of the set $\Lambda$. More accurate estimate of the error set in future algorithms may bring significant benefits for GSSC.
\par
As we mentioned,  erasures recovery is easier than error correction.
Putting new entries into the map of errors, we, in fact, announce them erasures. If we really have solid justification for this action, then the new iteration of GSSC can be considered as a new problem which is easier for SSC than the previous iteration.
\par
We will come back to the discussion about  interconnection between 
erasures and errors after presentation of numerical experiments. 
\par

\section{Numerical Experiments} 
\label{exper}
We will present the comparison of GSSC and SSC on synthetic data.
The input data was composed in accordance with the model given in \cite{EV}. 105 data vectors of dimension $D=50$ are equally split between  three 
4-dimensional linear spaces $\{\cS^i\}_{i=1}^3$. To make the problem more complicated each of those 3 spaces belongs to sum of two others. The smallest angles between spaces $\cS^i$ and $\cS^j$ 
are defined by formulas 
$$
\cos \theta_{ij}=\max_{\bu\in\cS^i,\,\bv\in\cS^j} \frac{\bu^T\bv}{\|\bu\|_2 \|\bv\|_2},\,i,j=1,2,3.
$$
We construct the data sets using vectors generated by decompositions with random coefficients in orthonormal bases $\be^j_i$ of spaces $\cS^j$. Three vectors $\be^j_1$ belongs to the same  2D-plane with angles $\widehat{\be^1_1 \be_ 1^2}=\widehat{\be^2_1 \be_ 1^3}=\theta$ and 
$\widehat{\be^1_1 \be_ 1^3}=2\theta$. The vectors $\be_2^1, \be_2^2, \be_3^1, \be_3^2, \be_4^1, \be_4^2$ are 
mutually orthogonal and orthogonal to $\be_1^1, \be_1^2, \be_1^3$; 
$\be_j^3=(\be_j^1+\be_j^2)/\sqrt{2}$, $j=2,3,4$.
The generator of standard normal distribution is used to generate data decomposition coefficients.
After the  generation, a random unitary matrix is applied to the result to avoid zeros in some regions of the matrix $Y$. 
\par
We use the notation $P_{ers}$  and $P_{err}$ for probabilities of erasures and  errors correspondingly.
\par
When we generate erasures we set random entries of the matrix $Y$ with probability $P_{ers}$ to zero and set those elements of $\Lambda$ to $\kappa=10^{-4}$. 
\par
The coordinates of samples with errors are generated randomly with probability $P_{err}$.
We use the additive model of errors,  adding  values of errors to the correct entries of $Y$.
The magnitudes of errors are taken from  standard normal distribution. 
\par
We run 20 trials of GSSC and SSC algorithms  for each combination of $(\theta, P_{err}, P_{ers})$, $$0\le\theta\le 60^\circ,$$ 
$$0\le P_{err}\le0.26,$$  
$$0\le P_{ers}\le0.4,$$ 
and output average values of misclassification. We note that for the angle $\theta=0$ the spaces $\{\cS^l\}$ have a common line and $\dim (\oplus_{l=1}^3\cS^l)=7$. Nevertheless, we will see that SSC and especially GSSC shows high capability  even for this hard settings.
\par
Now we describe the algorithm parameters. We do not use any creative stop criterion for greedy iterations. We set just make 5 iterations in each of 20 trials for all combinations 
$(\theta, P_{err}, P_{ers})$. The set $\Lambda$  is updated after each iteration as described above.
\par
The parameters for greedy envelop loop are: $\alpha_1=0.4$, $\alpha_2=0.5$, $\beta=0.65$, 
\par
The input parameters of the basic SSC block are as follows. We set $\alpha_e=5$, $\alpha_z=50$, $\rho^0=10$, $\mu=1.05$, $\epsilon=0.001$. 
The results presented on Fig. 1 %\ref{Fig:nonoise}
 confirms that GSSC has much higher error resilience than SSC. For all models of input data and for both algorithms "the phase transition 
curve" is in fact straight line with the slope about 0.4. In particular, this means that the influence on clustering  of one error is approximately corresponds to the influence of 2.5 erasures. We can see that GSSC gives reliable  clustering when $P_{err}+0.4P_{ers} \le 0.17$ for $\theta=60^\circ$ and $P_{err}+0.4P_{ers} \le 0.12$  for $\theta=6^\circ$.
For the case $\theta=0^\circ$, clustering cannot be absolutely perfect even for GSSC.  Indeed, the clustering for $P_{err}>0$ and $P_{ers}>0$ cannot be better than for error free model. At the same time, GSSC was designed for better error handling. In the error free case, GSSC has no advantage over SSC algorithm.  Thus, the images on Fig.~1 for $\theta=0$ have the gray background of the approximate level $0.042$ equal to the rate of misclassification of SSC. We believe that the reason of the misclassification lies in the method how we define the success. For $\theta=0^\circ$, there is a~common line belonging to all spaces $\cS^l$. For  points close to that line,  the considered algorithm has to make a hard decision, appointing only one cluster for each such point. 
\par
\noindent
\begin{picture}(00,110)
%\label{Fig2}
\put(0,5){          \epsfig{file=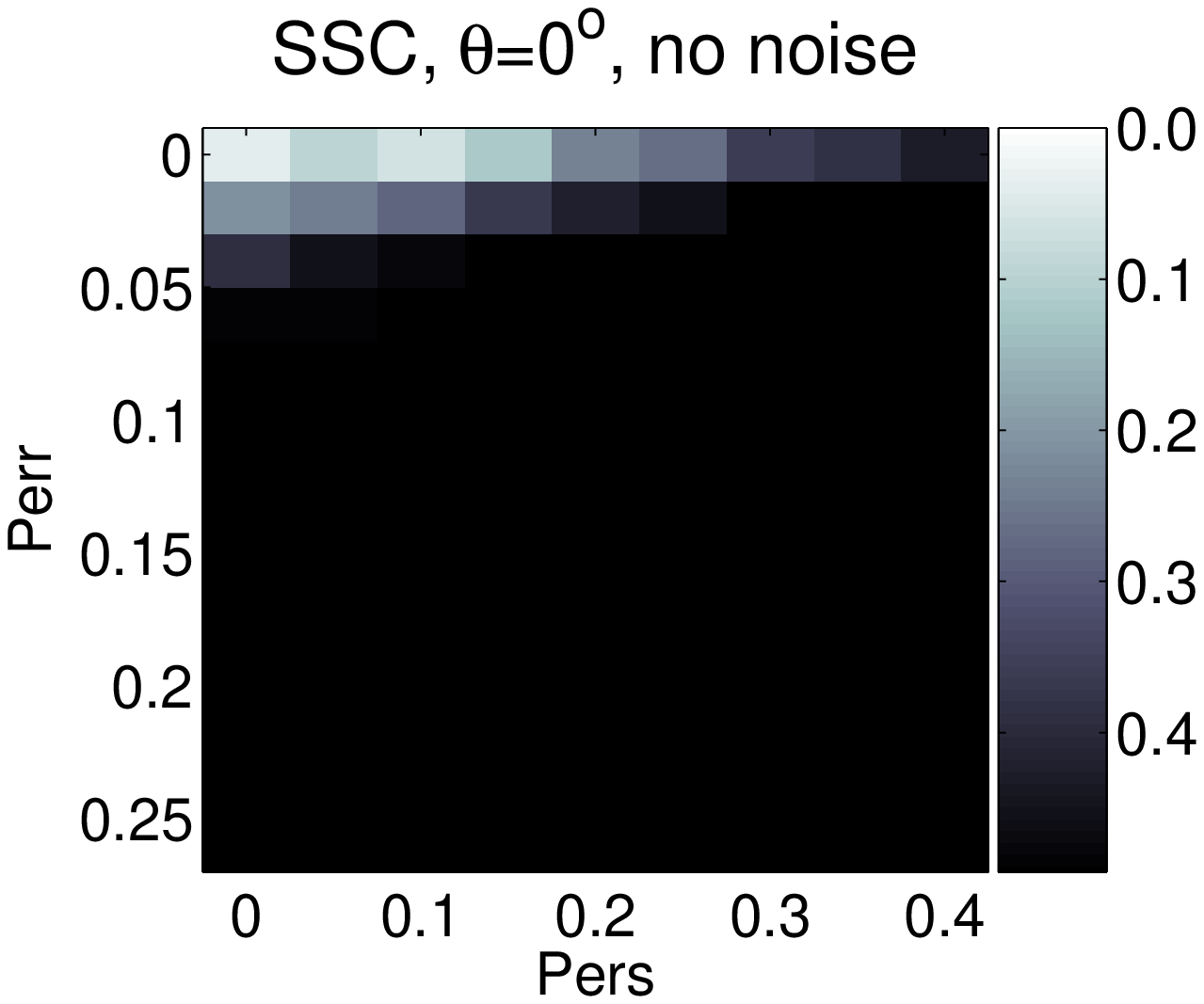,
           height=100pt, width=130pt}
		  } 
\put(120,5){          \epsfig{file=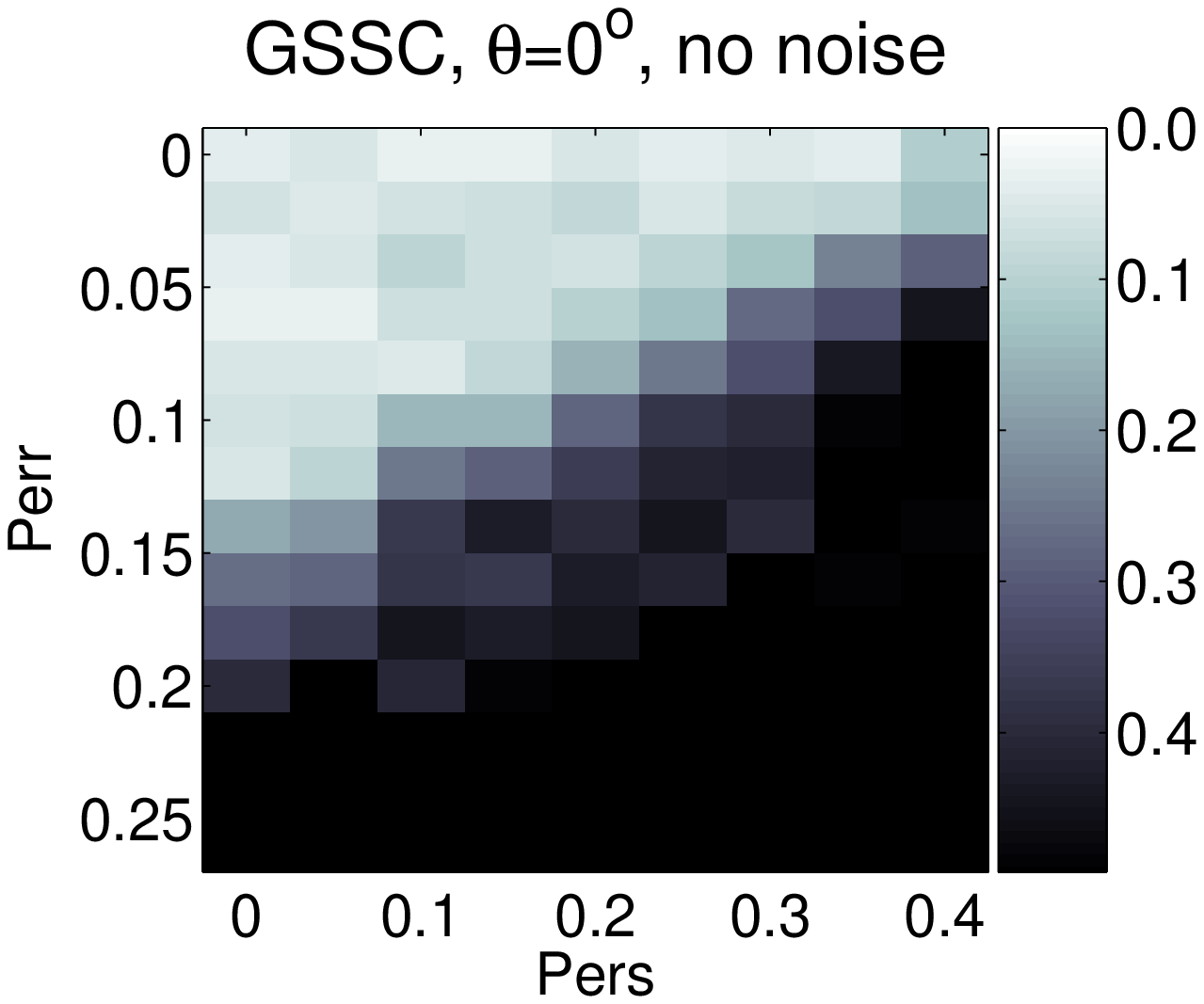,
           height=100pt, width=130pt}
		  } 
\end{picture}
\par
\noindent
\begin{picture}(00,110)
%\label{Fig2}
\put(0,5){          \epsfig{file=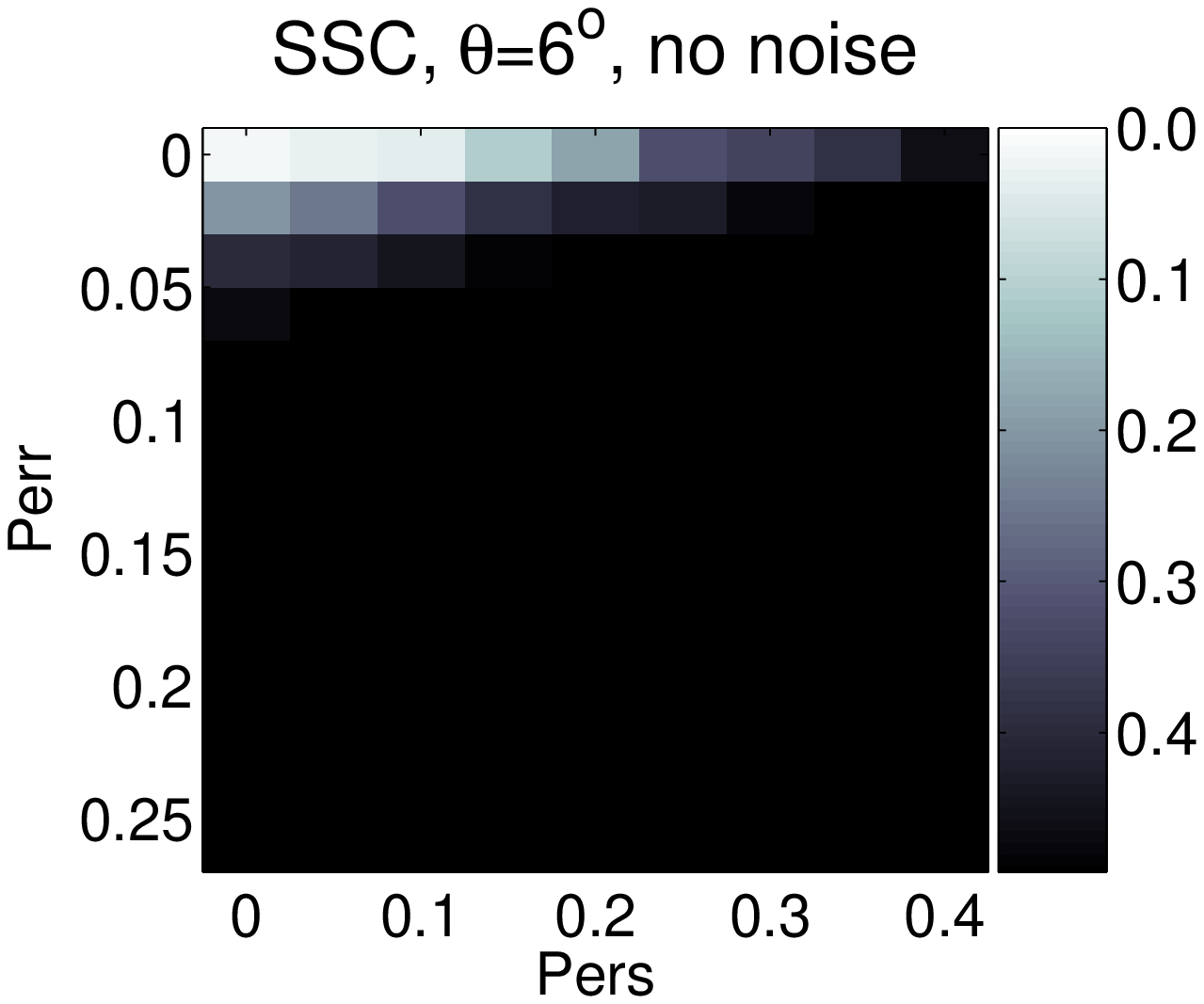,
           height=100pt, width=130pt}
		  } 
\put(120,5){          \epsfig{file=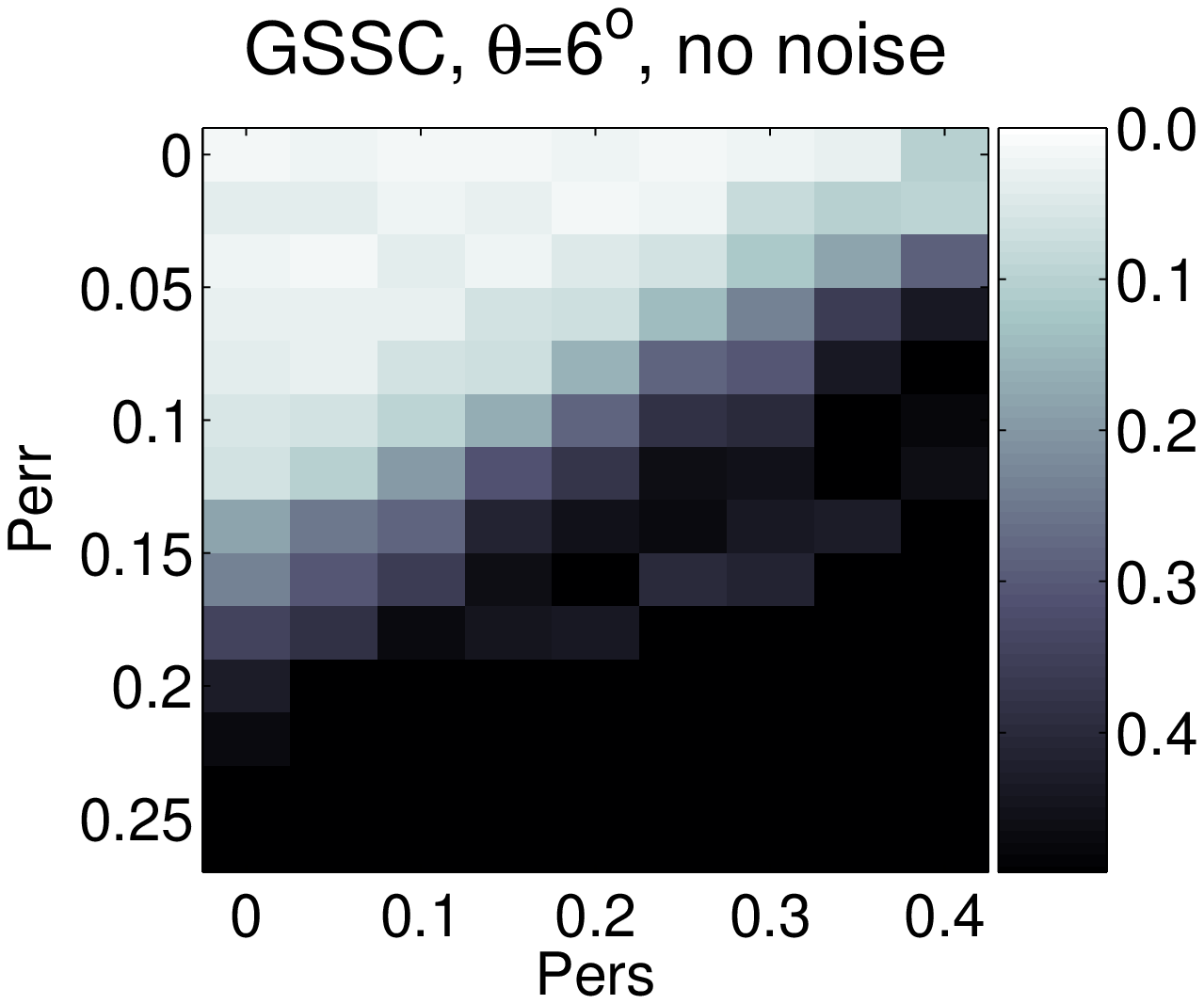,
           height=100pt, width=130pt}
		  } 
\end{picture}
\par
\noindent
\begin{picture}(00,115)
%\label{Fig2}
\put(0,10){          \epsfig{file=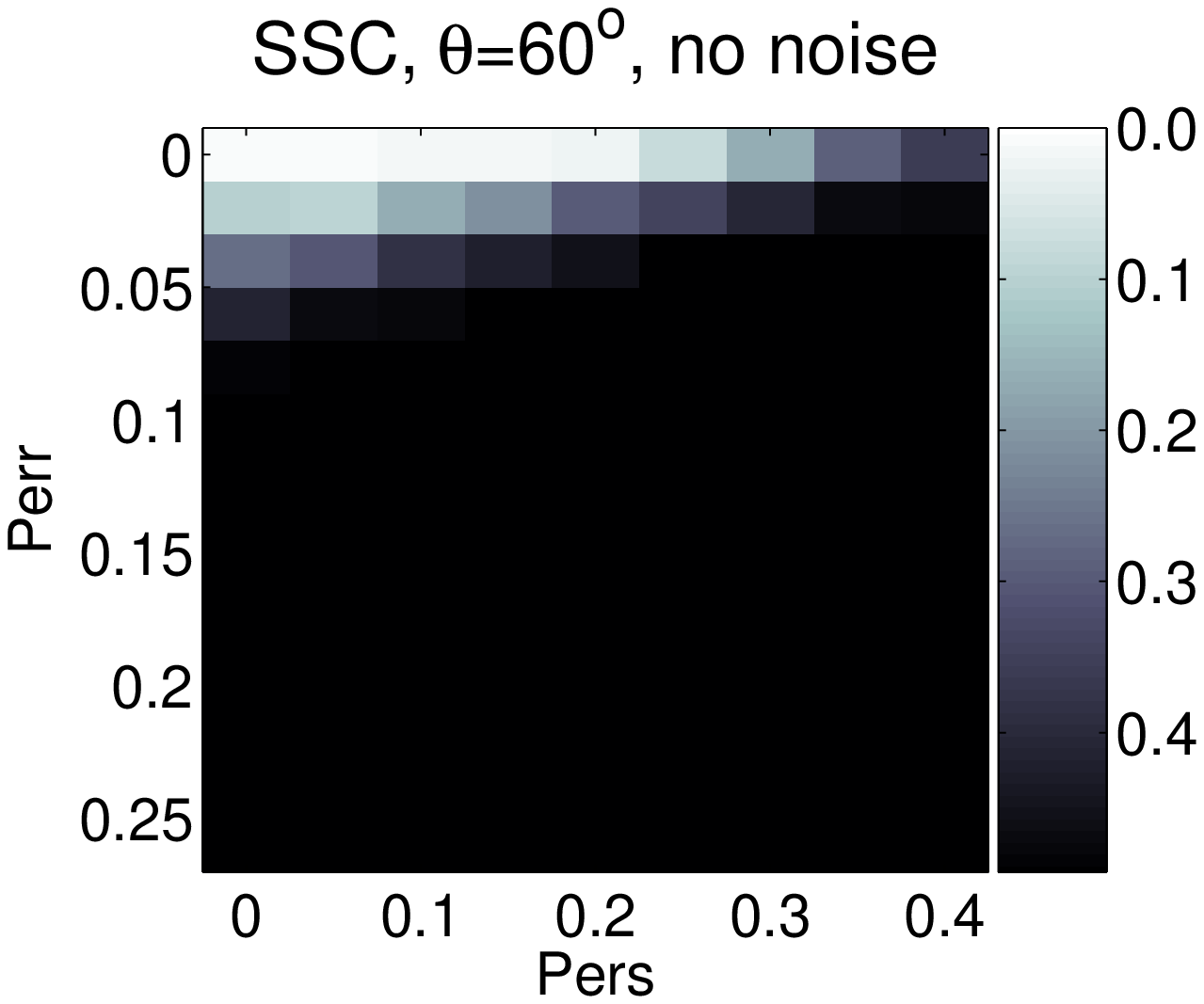,
           height=100pt, width=130pt}
		  } 
\put(120,10){          \epsfig{file=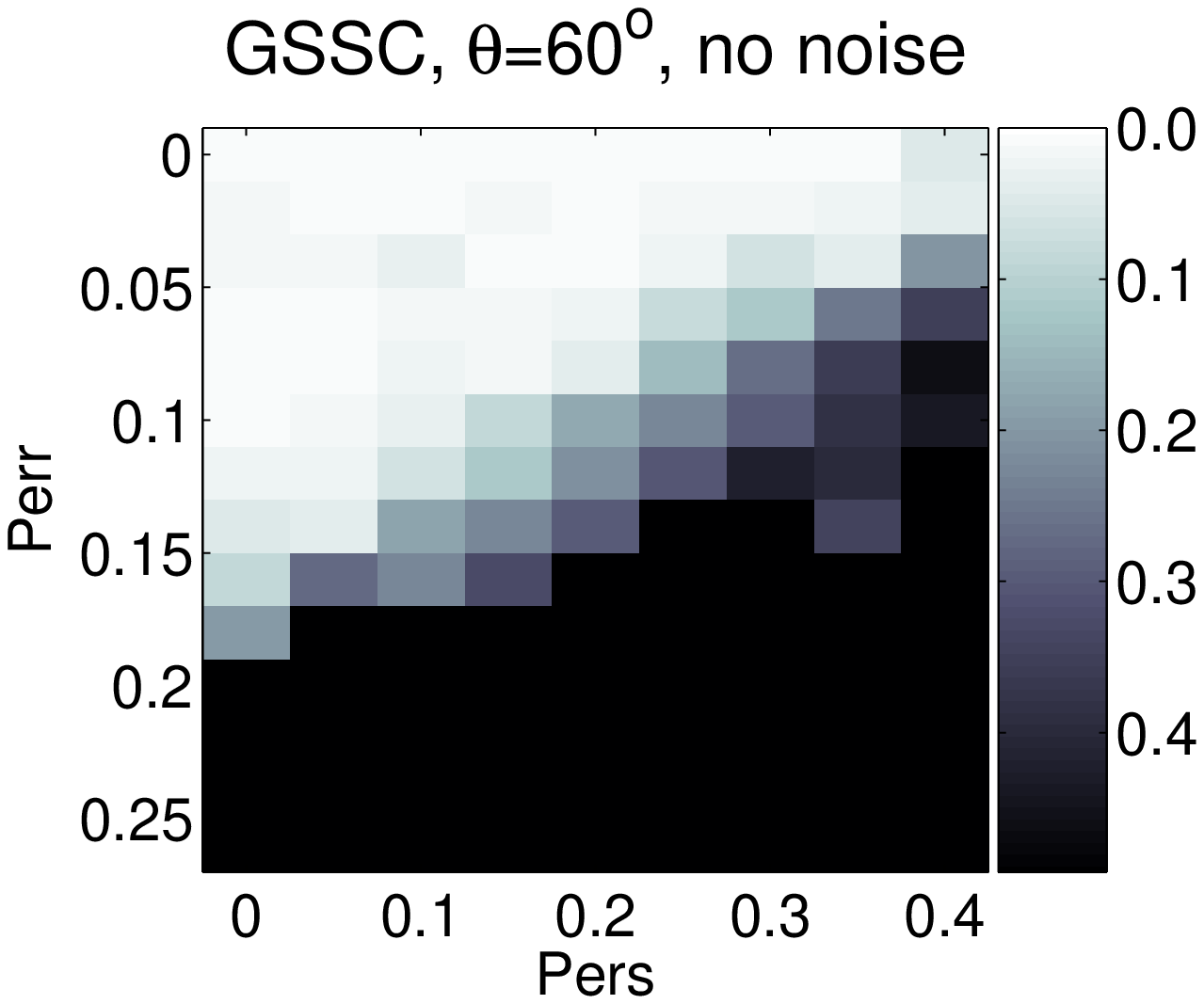,
           height=100pt, width=130pt}
		  } 
\put(70, 0){\footnotesize{\bf Fig.1. Clustering on noise free input.} }
%\label{Fig:nonoise}
\end{picture}
\par\noindent
Probably, for most of applied problems, the information about multiple point accommodation is more useful than unique space selection.
We advocate for such multiple selection because the typical follow-up problem after clustering is 
correction of errors in each of clusters. For this problem, it is not important to which of clusters the vector belonged from the beginning. When the vector affiliation is really important, side information has to be attracted.
\par
The second part of experiment deals with noisy data processing.  We apply to the matrix $Y$  independent Gaussian noise
of magnitude 10\% of mean square value of the data matrix $Y$, i.e., the noise level is -20dB. On Fig. 2, we present the results of processing of the noisy input analogous to results on Fig. 1. Evidently, that this quite strong noise has minor influence on the clustering efficiency. 
\par
If we increase the  noise up to -15 dB, the algorithms still resist. For -10 dB (see Fig. 3) GSSC  looses a  lot but still outperforms SSC.  Those losses are obviously caused by the increase of the noise fraction  in the mixture errors-erasures-noise, while the greedy idea efficiently works for highly localized corruption like errors and erasures.
\par
We emphasize that all results on Figs. 1--3 were obtained  with the same algorithm parameters. 
\par
\noindent
\begin{picture}(00,110)
%\label{Fig2}
\put(0,5){          \epsfig{file=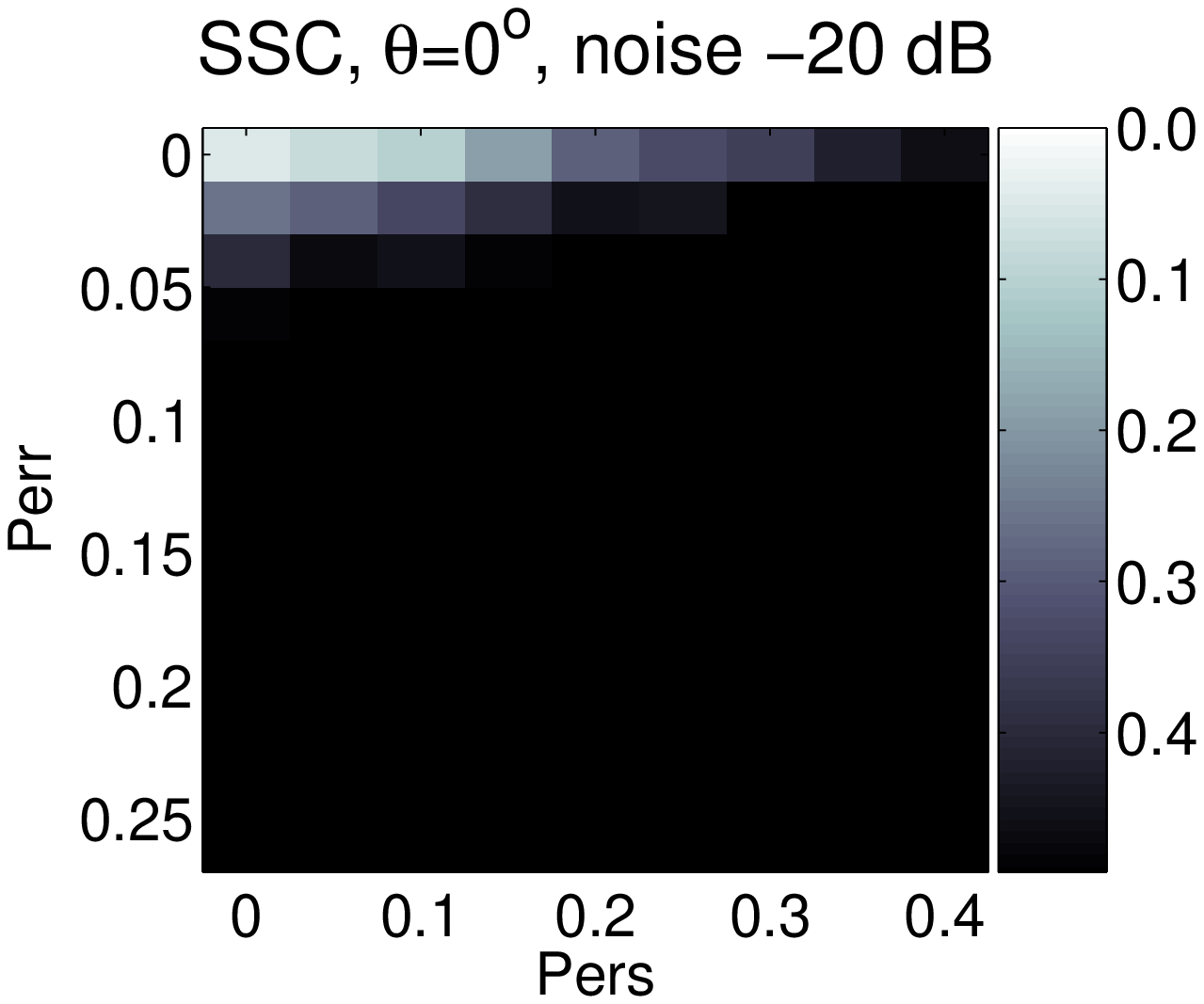,
           height=100pt, width=130pt}
		  } 
\put(120,5){          \epsfig{file=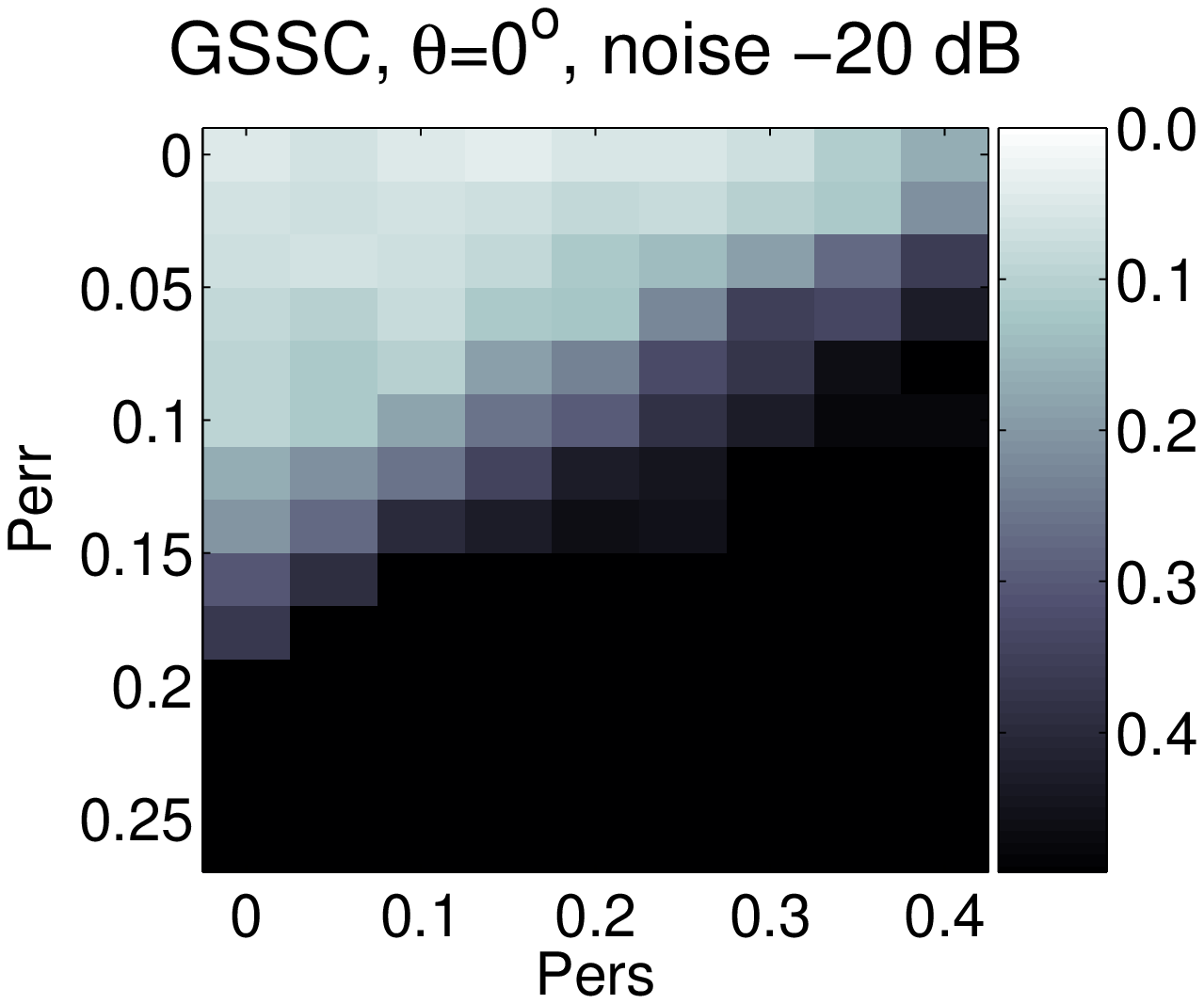,
           height=100pt, width=130pt}
		  } 
%\put(50, 5){\footnotesize {\bf Fig.2. Admissible error rate for erasure rate 0.1.}}
\end{picture}
\par
\noindent
\begin{picture}(00,110)
%\label{Fig2}
\put(0,5){          \epsfig{file=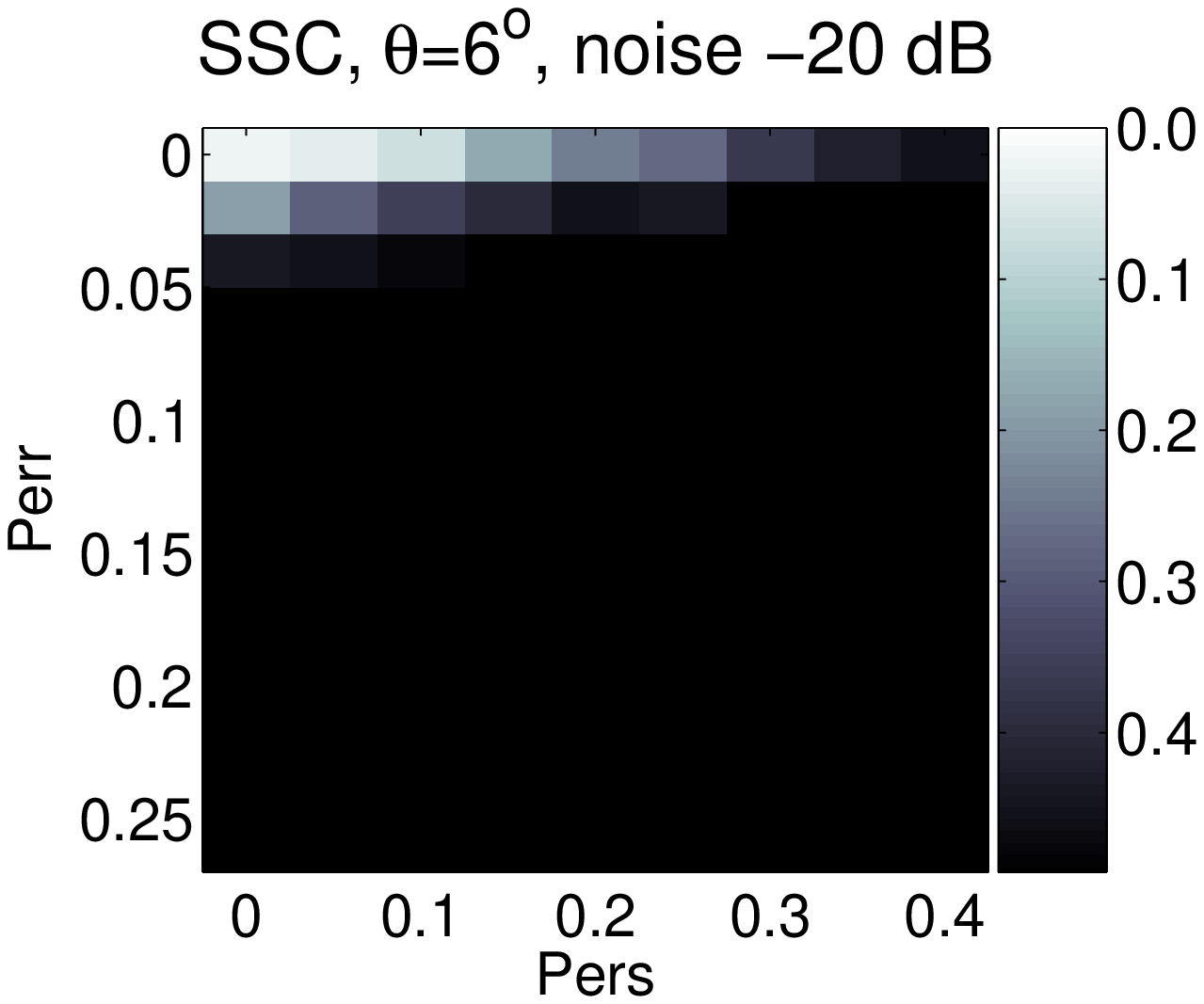,
           height=100pt, width=130pt}
		  } 
\put(120,5){          \epsfig{file=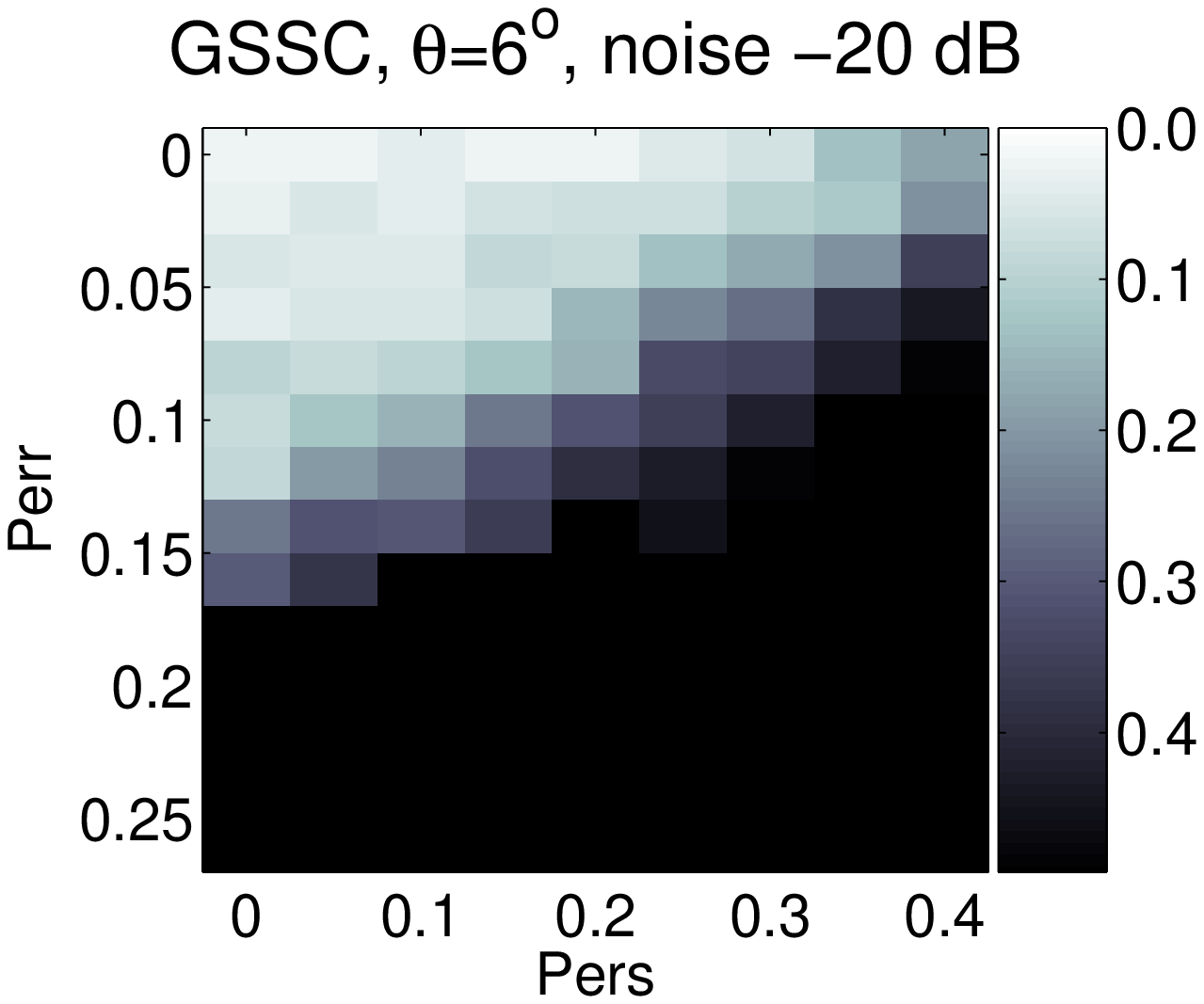,
           height=100pt, width=130pt}
		  } 
%\label{Fig:nonoise}
%\put(50, 5){\footnotesize {\bf Fig.2. Admissible error rate for erasure rate 0.1.}}
\end{picture}
\par
\noindent
\begin{picture}(0,115)
%\label{Fig2}
\put(0,10){          \epsfig{file=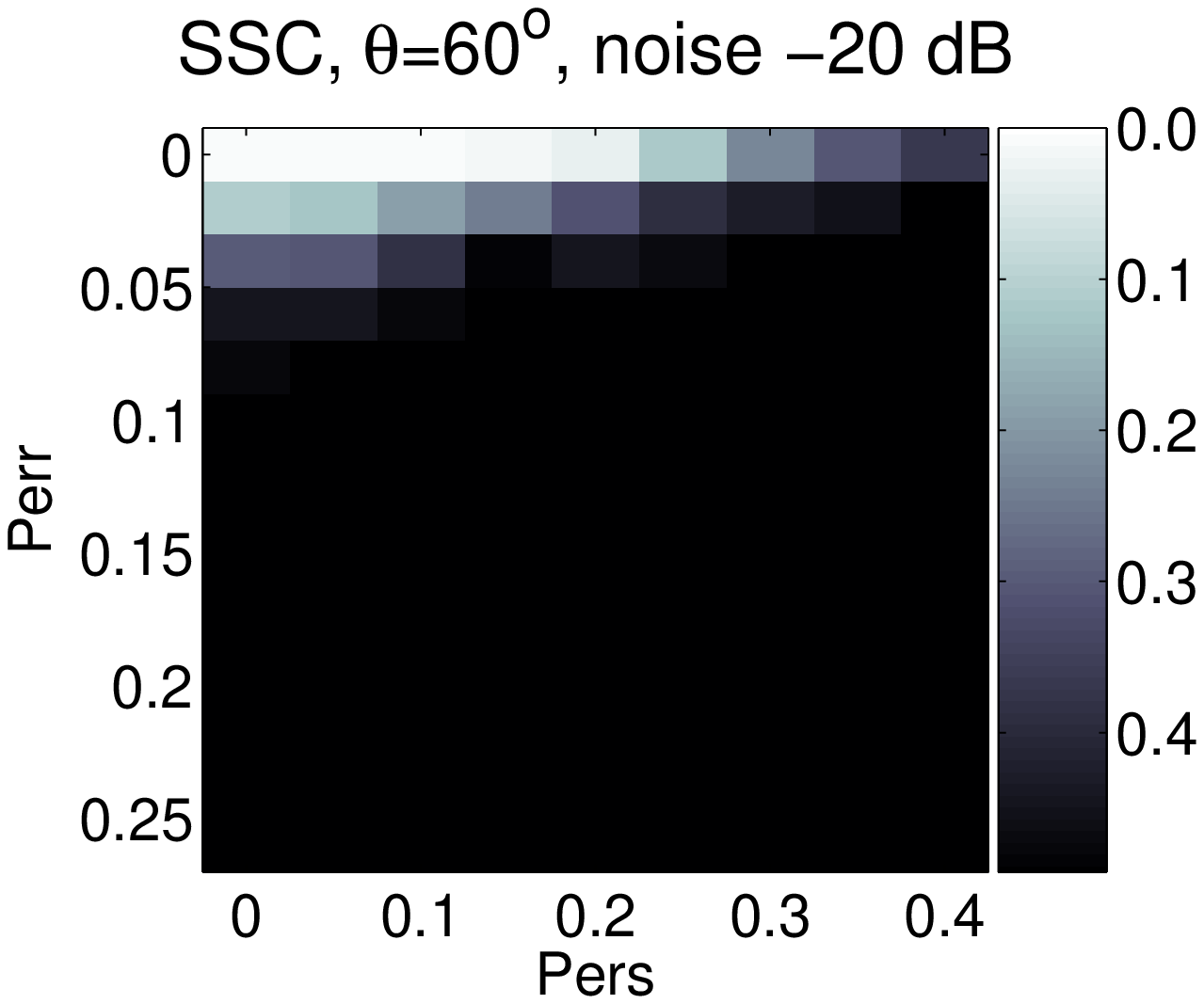,
           height=100pt, width=130pt}
		  } 
\put(120,10){          \epsfig{file=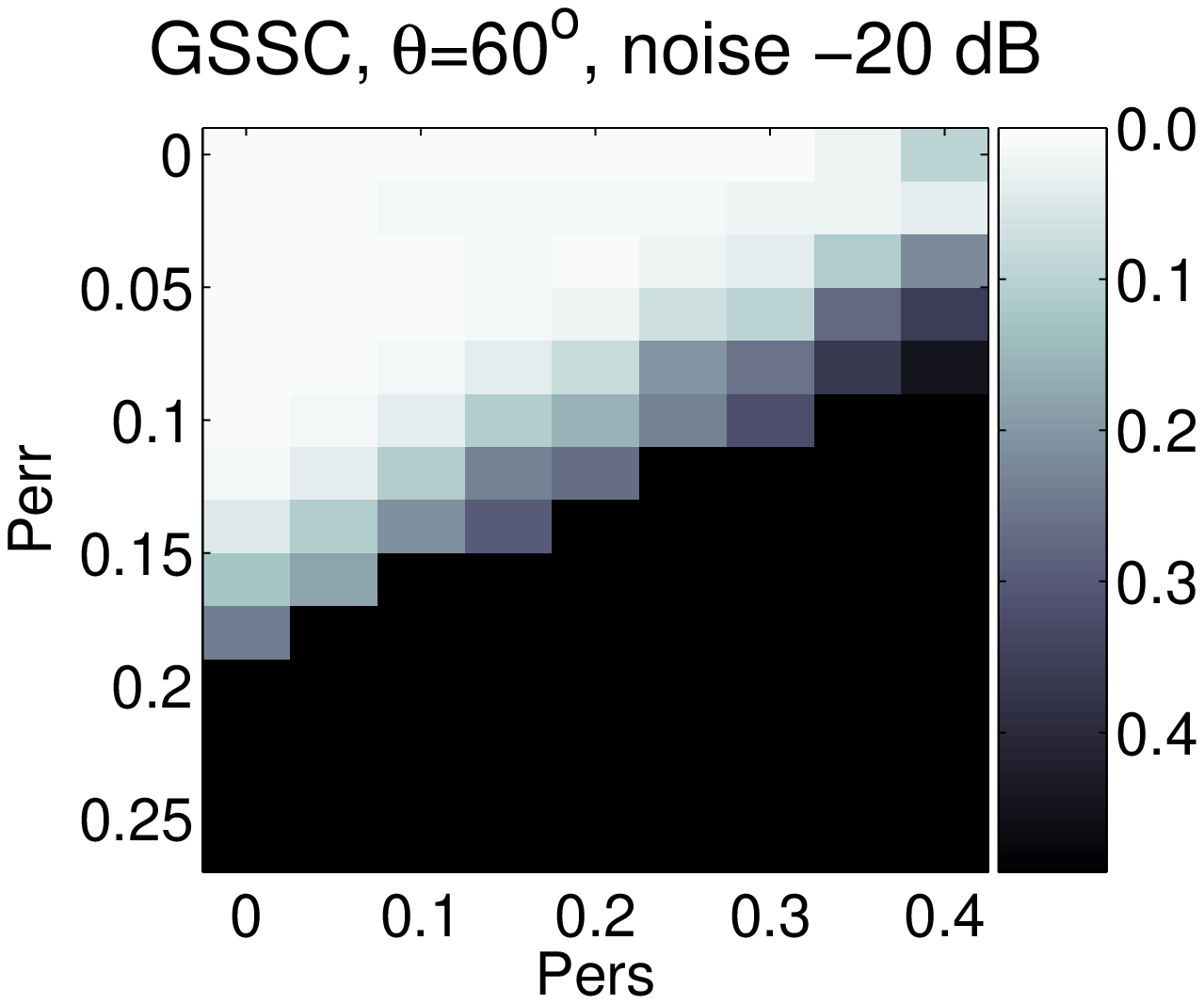,
           height=100pt, width=130pt}
		  } 
\put(50, 0){\footnotesize{\bf Fig.2. Clustering on noisy input, SNR=20 dB.} }
%\label{Fig:nonoise}
\end{picture}

\par
\noindent
\begin{picture}(00,115)
%\label{Fig2}
\put(0,10){          \epsfig{file=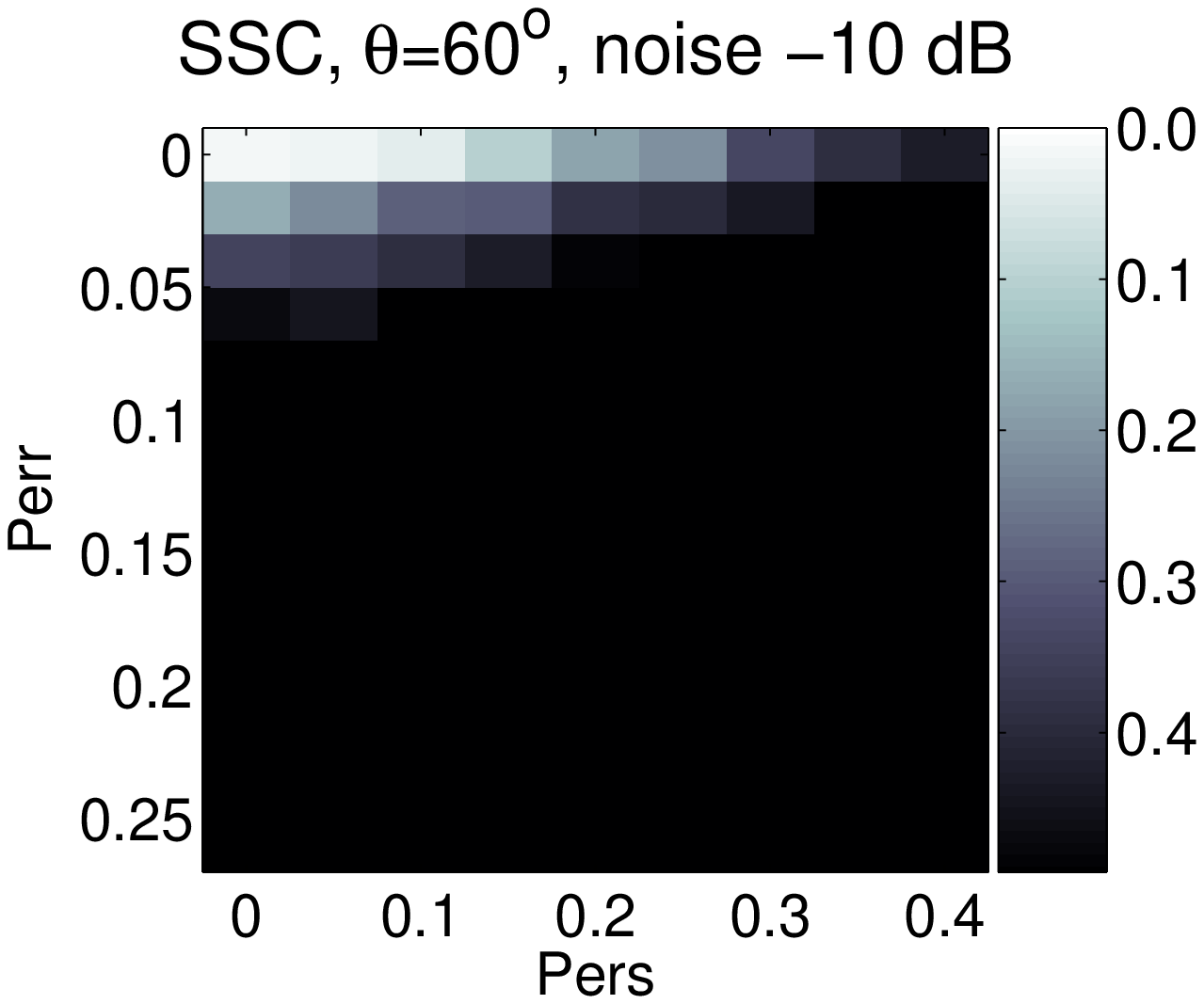,
           height=100pt, width=130pt}
		  } 
\put(120,10){          \epsfig{file=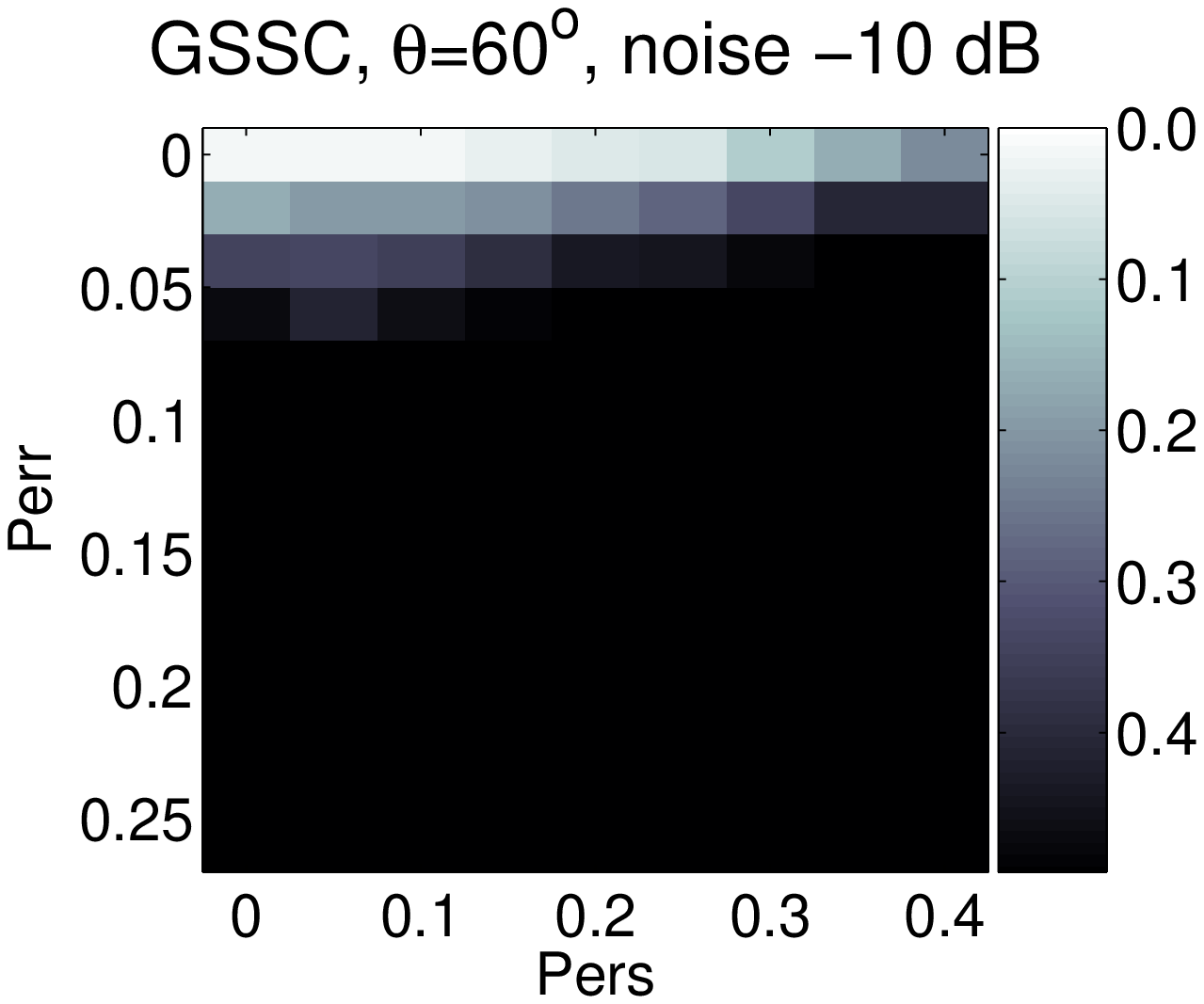,
           height=100pt, width=130pt}
		  } 
\put(50, 0){\footnotesize{\bf Fig.3. Clustering on noisy input, SNR=10 dB.} }
%\label{Fig:nonoise}
\end{picture}
\par
In conclusion, we demonstrate the role and dynamic of the greedy iterations.
In Tab. \ref{tab}, we give the dependence of the rate of misclassification for $P_{err}=0.05$, $P_{ers}=0.15$, SNR=20 dB from the number of GSSC iterations. The value Iter=0 corresponds to the pure SSC algorithm. The result was obtained as the average value of 100 trials.

\begin{table}[h]
\caption{Misclassification Rate  and Number of Iterations } 
\centering % centering table
\begin{tabular}{cccccccc} % creating eight columns
\hline
Iter  &0 &1&2&3&4&5&6\\
\hline 
$\theta=0^\circ$ &0.510 &0.479 &0.351  &0.199  &0.106 &0.053 &0.067\\
$\theta=60^\circ$ &0.465 &0.389  &0.237 &0.053 &0.022  &0.012 &0.010\\
\hline 
\end{tabular}
\label{tab}
\end{table}

\section{Conclusions and Future Work}
We consider a modification of Sparse Subspace Clustering algorithm based on greedy approach giving significant improvement of the algorithm resilience to
(simultaneous) entry corruption, incompleteness, and noise. 
\par
While the basic SSC algorithm has some internal resilience to corruption, it reduces error and erasure influence on clustering quality, strictly speaking, it does not have error correction capabilities. 
\par
Adding the error correction capability  may have independent importance as well as it may improve the clustering quality. This direction deserves the further research.
\par
In the described version, GSSC  has 5-6 iterations of the algorithm SSC, having proportional increase of computing time. We believe that this computing time increase can be significantly eliminated with preserving the algorithm efficiency if updates of $\Lambda$ are incorporated into internal iterations of SSC algorithm. This option was implemented by us in very analogous situation for acceleration of the $\ell^1$-greedy algorithm in \cite{PK2}.  Finding an appropriate stopping  criterion also could 
reduce computing time and improve clustering.
\par
One more reserve for algorithm improvement is selection of the parameters adaptive to input data.
While, as we mentioned all results were obtained with the same set of parameters, the adaptation may bring significant increase of algorithm capability. One of such adaptive solution for error correction in Compressed  Sensing was recently found by the authors in~\cite{PK3}.
%%%%%%%%%%%%%%%%%%%%%%%%%%%%%%%%%%%%%%%%


\begin{thebibliography}{99}
\providecommand{\url}[1]{#1}
\csname url@rmstyle\endcsname
\providecommand{\newblock}{\relaix}
\providecommand{\bibinfo}[2]{#2}
\providecommand\BIBentrySTDinterwordspacing{\spaceskip=0pt\relax}
\providecommand\BIBentryALTinterwordstretchfactor{4}
\providecommand\BIBentryALTinterwordspacing{\spaceskip=\fontdimen2\font plus
\BIBentryALTinterwordstretchfactor\fontdimen3\font minus
  \fontdimen4\font\relax}
\providecommand\BIBforeignlanguage[2]{{%
\expandafter\ifx\csname l@#1\endcsname\relax
\typeout{** WARNING: IEEEtran.bst: No hyphenation pattern has been}%
\typeout{** loaded for the language `#1'. Using the pattern for}%
\typeout{** the default language instead.}%
\else
\language=\csname l@#1\endcsname
\fi
#2}}


\bibitem{BPCPE} S.Boyd, N.Parikh, E.Chu, B.Peleato, and J. Eckstein, Distributed Optimization and Statistical Learning via the Alternating Direction Method of Multipliers,  Foundations and Trends in Machine Learning, V.3 (2010), No 1,  1--122.


\bibitem{CR}
E.J.~Cand\`es  and B.~Recht,
Exact Matrix Completion via Convex Optimization, Comm. of the ACM, V. 55, \# 6,  2012, 
111--119 

\bibitem{CT}
E.J. Cand\`es, T. Tao, Decoding by linear programming, IEEE Transactions on Information
Theory, 51 (2005), 4203--4215.


\bibitem{CWB}
E. J. Candes, M. B. Wakin, and S. Boyd,
 \emph{Enhancing sparsity by reweighted $\ell_1$ minimization}
J. of Fourier Anal. and Appl., special issue on sparsity, {\bf 14} (2008), 877--905 .


\bibitem{CJSC}
Y. Chen, A. Jalali, S. Sanghavi and C. Caramanis, Low-rank matrix recovery from errors and erasures, IEEE International Symposium on Information Theory, Proceedings (ISIT), 2011,  2313--2317 


\bibitem{D}
D. Donoho, Compressed Sensing, IEEE Trans. on Information Theory, 52 (2006), 1289--1306.



\bibitem{EV0}
E.Elhamifar, R.Vidal, Sparse  Subspace Clustering, IEEE Conference on Computer Vision and Pattern Recognition,  20-25 June 2009, 2790--2797 

\bibitem{EV}
E.Elhamifar, R.Vidal, Sparse  Subspace Clustering: Algorithm, Theory, and Applications, arXiv:2013.1005v3 [cs.CV], 5 Feb. 2013.

\bibitem{KP1}
 I.Kozlov, A.Petukhov,
Sparse Solutions for Underdetermined Systems of Linear Equations, chapter in “Handbook of Geomathematics”, Springer, 1243--1259, 2010.

\bibitem{LCM}
Z. Lin, M. Chen, Yi Ma
The Augmented Lagrange Multiplier Method for Exact Recovery of Corrupted Low-Rank Matrices, Preprint, arXiv:1009.5055, 2010; rev. 9 Mar 2011.

\bibitem{L} U. von Luxburg, A Tutorial on Spectral Clustering, Statistics and Computing, 17 (2007), 25p.


\bibitem{NWJ} A.Ng, Y.Weiss, M.Jordan, On Spectral Clustering: Analysis and an Algorithm, Neural Information Processing Sustems, 2001, 849--856.

\bibitem{PK2}
A.Petukhov, I.Kozlov, 
Fast Implementation of $\ell^1$-greedy algorithm, Recent Advances in Harmonic Analysis and Applications, Springer, 2012, 317--326.

\bibitem{PK3}
A.Petukhov, I.Kozlov, Correcting Errors in Linear Measurements and Compressed Sensing of Multiple Sources, submitted to Applied Mathematics and Computation, 2012.

\bibitem{PK4}
A.Petukhov, I.Kozlov,  Greedy Approach for Low-Rank matrix recovery, Submitted to WorldComp 2013, Proceedings.

\bibitem{RTVM}
S.Rao, R.Tron, R.Vidal, Yi Ma, Motion Segmentation in the Presence of Outlying, Incomplete, or Corrupted Trajectories, 
IEEE Transactions on Pattern Analysis and Machine Intelligence, 32(2010), 1832--1845.

\bibitem{RV}
M.Rudelson, R.Vershynin, Geometric approach to error correcting codes and reconstruction
of signals, International Mathematical Research Notices 64 (2005), 4019--4041.

\bibitem{WYGSM}
J.Wright, A. Y. Yang, A. Ganesh, S. S. Sastry, and Yi Ma, Robust Face Recognition
via Sparse Representation, Pattern Analysis and Machine Intelligence, IEEE Transactions on,
Feb 2009, Volume: 31 , Issue: 2 Page(s): 210--227.

\bibitem{YYO}
M.Yan, Y.Yang, S.Osher, Exact Low-Rank Matrix Completion from Sparsely Corrupted Entries Via Adaptive Outlier Pursuit, J. Sci. Comput., January 2013.

\end{thebibliography}
\end{document}